\documentclass{article}
\usepackage{amsmath}
\usepackage{amsfonts}

\newenvironment{proof}[1][Proof]{\noindent\textbf{#1.} }{\ \rule{0.5em}{0.5em}}
\input{tcilatex}
\begin{document}

\begin{center}
\begin{equation*}
\text{{\LARGE \ }On Salem numbers \ which are exceptional units}
\end{equation*}

\bigskip\ Toufik Za\"{\i}mi
\end{center}

\textbf{Abstract. }\textit{By extending a construction due to \ Gross and
McMullen }[6],\textit{\ we show that for any odd integer }$n$ \textit{and
for any even integer} $d\geq n+3$ \textit{there are infinitely many Salem
numbers }$\alpha $\textit{\ of degree }$d$\textit{\ such that }$\alpha ^{n}-1
$\textit{\ is a unit. A similar result is proved when} $n$\textit{\ runs
through some classes of even integers, }$d\geq n+4$ \textit{and }$d/2$ 
\textit{is an odd integer.}

\medskip



\smallskip

\begin{center}
\textbf{1. Introduction}
\end{center}

A Salem number is a real algebraic integer greater than $1,$ whose other
conjugates lie inside the closed unit disc, with at least one conjugate
lying on the boundary; the set of such numbers is traditionally denoted by $%
\mathbb{T\ }$[2].

Let $S_{\alpha }=S$ \ designate the minimal polynomial\ of a Salem number $%
\alpha .$ Then, $S$ has two real roots, namely $\alpha _{1}:=\alpha $ and $%
\alpha ^{-1},$ and the rest, say $\alpha _{2}^{\pm 1},...,\alpha _{t}^{\pm
1},$ lie on the unit circle. Also, $\deg (S)=2t\geq 4,$ $S$ is reciprocal,
i. e., $S(x)=x^{t}S(1/x),$ and the polynomial $T_{\alpha
}(x)=T(x):=(x-(\alpha _{1}+\alpha _{1}^{-1}))(x-(\alpha _{2}+\overline{%
\alpha _{2}}))\cdot \cdot \cdot (x-(\alpha _{t}+\overline{\alpha _{t}}))\in 
\mathbb{Z}[x],$ called the trace polynomial of $S$ (or the trace polynomial
of $\alpha )$ satisfies 
\begin{equation}
S(x)=x^{t}T(x+\frac{1}{x}).  \tag{1}
\end{equation}%
Clearly, $T_{\alpha }$ \ is irreducible and has one root greater than $2,$
namely $\alpha +\alpha ^{-1}$ (called in [6] a Salem trace number), and $%
(t-1)\geq 1$ roots belonging to the interval $(-2,2).$ Conversely, an
algebraic integer $\beta >2$ of degree $t\geq 2$ whose other conjugates lie
in $(-2,2)$ is a Salem trace number associated to a Salem number $\alpha $
of degree $2t,$ via the relation $\beta =\alpha +\alpha ^{-1}.$ Some
properties of Salem trace polynomials may be found in [10].

An algebraic integer $u$ is said to be an exceptional unit if both $u$ and $%
u-1$ are units [7]. Siegel proved that there are at most a finite number of
exceptional units in any number field $F,$ and a result of Evertse [5]
implies that this number is no more than $3(7^{3\deg (F)}).$

Clearly, a Salem number is a unit. Theorem 3 of [8] yields that every
integer greater than $2$ is a limit of a sequence of Salem numbers which are
exceptional units, and a characterization of Salem numbers $\alpha $ that
are exceptional units among the elements of \ the set $\mathbb{T},$ in terms
of the fractional parts of the powers of $\alpha ,$ follows from the main
result of [9].

In [7], Silverman conducted some investigations suggesting that if a unit $u$
has small Mahler measure (the Mahler measure of an algebraic integer is the
absolute value of the product of its conjugates of modulus greater than $1),$
then there are many values of $n$ for which $u^{n}-1$ is a unit. In
particular, he proved that the numbers of such values of $n$ is bounded \
above by $O(d^{1+0.7/\log \log d}),$ where $d=\deg (u).$ In most of his
numerical investigations Silverman supposed that $u$ is a Salem number from
Boyd's lists [3-4], and one of his observations \ says that among the forty
two known small Salem numbers thirty three are exceptional units, and there
is a Salem number $\alpha $ of degree $26$ such that $\alpha ^{n}-1$ is a
unit for all $n\in \{1,2,...,10\}.$ From Theorem 2, below, we can deduce
that for any finite set $F$ of natural numbers there is an infinite subset $N
$ of $2\mathbb{N}$ such that for each $d\in N$ \ there exist infinitely many
Salem numbers $\alpha $ of degree $d$ having the property that $\alpha ^{n}-1
$ is a unit for all $n\in F.$

By considering the problem of the realization of a monic irreducible
polynomial with integer coefficients as the characteristic polynomial of an
automorphism of the even indefinite unimodular lattice, Gross and McMullen
[6] were concerned with the elements $\alpha $ of $\mathbb{T}$ such that $%
\alpha ^{2}-1$ is a unit. They called such numbers $\alpha $ unramified
Salem numbers, and they showed in this case that $\deg (\alpha )/2$ must be
odd [6, Proposition 3.3]. Also, they gave a construction leading to the fact
that for any odd integer $t\geq 3$ there are infinitely many unramified
Salem numbers of degree $2t.$

Since a power of a unit (resp. of a Salem number of degree $2t$) is a unit
(resp. is a Salem number of degree $2t$), a Salem number is unramified if
and only if its square is an exceptional unit, and a natural question arises
immediately:

\smallskip

\textit{For which integers }$n\geq 1$\textit{\ and }$t\geq 2$\textit{\ are
there infinitely many Salem numbers }$\alpha $\textit{\ of degree }$2t$%
\textit{\ such that }$\alpha ^{n}-1$\textit{\ is a unit?}

\smallskip

The aim of the present paper is to give some partial answers to this
question.

\medskip

\textbf{Theorem 1. }\textit{There exist infinitely many Salem numbers }$%
\alpha $\textit{\ of degree }$2t$\textit{\ such that }$\alpha ^{n}-1$\textit{%
\ is a unit, whenever one of the following conditions holds:}

\textit{(i) }$n$\textit{\ is odd and }$t\geq (n+3)/2;$

\textit{(ii) }$n\equiv 2\func{mod}4,$\textit{\ }$t$\textit{\ is odd and }$%
t\geq (n+4)/2;$

\textit{(iii) }$n=2^{s}$\textit{\ for some integer }$s\geq 2,$\textit{\ }$t$%
\textit{\ is odd and }$t\geq (n+6)/2;$

\textit{(iv) }$n\equiv 4\func{mod}8,$ $n\neq 0\func{mod}3,$\textit{\ }$t$%
\textit{\ is odd and }$t\geq (n+6)/2.$

\bigskip

Because the degree $t$ \ of the polynomial trace is at least $2,$ Theorem
1(i) is sharp when $n=1.$ According to [6, Proposition 3], Theorem 1(ii) is
also sharp for $n=2.$ As mentioned above, [6, Proposition 3.3] says that if $%
\alpha \in \mathbb{T}$ and $\alpha ^{n}-1$ is a unit for some even natural
natural $n,$ then $t$ is necessarily odd, since  $\alpha ^{n}-1=(\alpha
^{2}-1)((\alpha ^{2})^{n/2-1}+(\alpha ^{2})^{n/2-2}+\cdot \cdot \cdot +1).$
We will see in the next section that Theorem 1 is sharp when $n\in \{3,4,6\}$
but not for $n=5.$ The following result is a weak complement of Theorem 1
when $n\equiv 0\func{mod}4.$

\bigskip

\textbf{Theorem 2. }\textit{Let} $n\in \mathbb{N}.$ \textit{Then, there is
an infinite subset }$N$ \textit{of } $\mathbb{N}$ \textit{such that for each 
}$t\in N$\textit{\ \ there are infinitely many Salem numbers }$\alpha $%
\textit{\ satisfying }$\deg (\alpha )=2t$\textit{\ and }$\alpha ^{n}-1$%
\textit{\ is a unit.}

\bigskip

Given a finite sequence $n_{1}<n_{2}<\cdot \cdot \cdot <n_{s}$ of natural
numbers, if we apply Theorem 2 with $n=n_{1}n_{2}...n_{s}$ we see that there
is an infinite subset (say again) $N$ of $2\mathbb{N}$ such that for each $%
d\in N$ \ there exist infinitely many Salem numbers $\alpha $ of degree $d$
having the property that $\alpha ^{n_{k}}-1$ is a unit for each $k\in
\{1,...s\}.$

The first value of $n$ for which Theorem 1 does not answer our question is $%
12.$ From the proof of Theorem 2, it is easy to see that for any natural
number $v\neq 0\func{mod}3$ there are infinitely many Salem numbers $\alpha $
of degree $18+4v$ such that $\alpha ^{12}-1$ is a unit (the problem remains
open for each degree of the form $18+4v$ with $v\equiv 0\func{mod}3).$

Throughout, when we speak about conjugates, the norm, the degree and the
minimal polynomial of an algebraic integer without mentioning the basic
field, this is meant over $\mathbb{Q}.$ By convention, an empty sum is equal
to zero and an empty product is equal to one. Also, we denote by%
\begin{equation*}
C_{n}(x):=\dprod\limits_{k=1}^{(n-1)/2}(x-2\cos (\frac{2k\pi }{n}))\text{ \
\ \ (resp. }C_{n}(x):=\dprod\limits_{k=1}^{(n-2)/2}(x-2\cos (\frac{2k\pi }{n}%
)))
\end{equation*}%
the trace polynomial of $(x^{n}-1)/(x-1)$ when $n$ is odd (resp. of $%
(x^{n}-1)/(x^{2}-1)$ when $n$ is even). Then, $C_{1}(x)=C_{2}(x)=1,$ $%
C_{n}(x)\in \mathbb{Z}[x]$ and the roots of $C_{n}$ (for $n\geq 3)$ are
distinct and belong to the interval $(-2,2).$ To prove Theorem 1 and Theorem
2 we shall use the following generalization of [6, Theorem 7.3].

\bigskip

\ \textbf{Theorem 3. }\textit{Let }$n\geq 1$ \textit{be an odd }(\textit{%
resp. an even}) \textit{natural number, }$t$ \textit{an integer greater \
than or equal to} $(n+3)/2$ (\textit{resp.} $t$ \textit{an odd integer
greater than or equal to} $(n+4)/2)$ \textit{and } $D(x)\in \mathbb{Z}[x]$ 
\textit{a monic separable polynomial of degree }$t-(n+3)/2$ (\textit{resp.} 
\textit{of degree} $t-(n+4)/2)$ \textit{with all roots }(\textit{if any})%
\textit{\ in }$(-2,2)$ \textit{and such that \ }$\gcd (D,C_{n})=1.$ \textit{%
Then, for all sufficiently large} \textit{integers }$a$\textit{\ the
polynomial} 
\begin{equation*}
C_{n}(x)(x-2)D(x)(x-a)-1\text{ \ \ \ \ \ (\textit{resp.} }%
C_{n}(x)(x^{2}-4)D(x)(x-a)-1)
\end{equation*}%
\textit{is a trace polynomial }(\textit{of degree }$t)$ \textit{of a Salem
number }$\alpha $\textit{\ such that }$\alpha ^{n}-1$\textit{\ is a unit.}

\medskip \medskip\ 

The proofs of the theorems are given in the last section. In the next one we
consider Salem numbers $\alpha $ such that $\alpha ^{n}-1$ is a unit for
some $n\leq 6,$ and we give a characterization of Salem polynomials $%
S_{\alpha }$ such that $\alpha ^{n}-1$ is a unit and  $n\neq 5.$ All
computations are done using the system Pari [1].

\bigskip

\begin{center}
\textbf{2. On Salem numbers }$\alpha $ \textbf{such that }$\alpha ^{n}-1$%
\textbf{\ is a unit for some} $n\leq 6$
\end{center}

With the notation above, set 
\begin{equation*}
S(x)=(1+x^{2t})+a_{1}(x+x^{2t-1})+\cdot \cdot \cdot
+a_{t-1}(x^{t-1}+x^{t+1})+a_{t}x^{t},
\end{equation*}%
for the minimal polynomial of the Salem number $\alpha .$

Since the distinct numbers $\alpha _{1}^{\pm n},...,\alpha _{t}^{\pm n}$ $,$
where $n\in \mathbb{N},$ are the conjugates of the Salem number $\alpha
^{n}, $ we have $\func{Norm}(\alpha ^{n}-1)=(\alpha _{1}^{n}-1)(\alpha
_{1}^{-n}-1)(\alpha _{2}^{n}-1)(\overline{\alpha _{2}^{n}}-1)\cdot \cdot
\cdot (\alpha _{t}^{n}-1)(\overline{\alpha _{t}^{n}}-1)<0$ and $\func{Norm}%
(\alpha ^{n}+1)>0.$ Hence, $\alpha ^{n}-1$ is a unit if and only if $\func{%
Norm}(\alpha ^{n}-1)=-1,$ and $\alpha ^{n}+1$ is a unit if and only if $%
\func{Norm}(\alpha ^{n}+1)=1.$ Especially, 
\begin{equation}
\alpha -1\text{ is a unit }\Leftrightarrow S(1)=-1,  \tag{2}
\end{equation}%
as $S(1)=\func{Norm}(\alpha -1),$ or equivalently 
\begin{equation}
\alpha -1\text{ is a unit}\Leftrightarrow a_{t}=-3-2(a_{1}+\cdot \cdot \cdot
+a_{t-1}).  \tag{3}
\end{equation}%
In a similar way we obtain a characterization of the polynomial $S$ \ when $%
\alpha ^{n}-1$ is a unit for some $n\in \{2,3,4\}:$

\medskip

\medskip \textbf{Proposition 1.}\textit{\ With the notation above we have
the following assertions.}

\textit{(i)} $\alpha ^{2}-1$\textit{\ is a unit\ }\ \textit{if and only if} 
\textit{\ } 
\begin{equation}
t\equiv 1\func{mod}2,\text{ \ }a_{t}=-1-2\dsum\limits_{k=0}^{(t-3)/2}a_{2k+1}%
\text{ \textit{and }\ }a_{t-1}=-1-\dsum\limits_{k=1}^{(t-3)/2}a_{2k}; 
\tag{4}
\end{equation}

\textit{(ii) }$\alpha ^{4}-1$\textit{\ is a unit\ }\ \textit{if and only if }%
$t\equiv 1\func{mod}4,$\textit{\ } 
\begin{equation}
a_{t}=-1-2\dsum\limits_{k=0}^{(t-5)/4}a_{4k+1},\text{ }a_{t-1}=-1-\dsum%
\limits_{k=1}^{(t-3)/2}a_{2k}\text{ \textit{and }}a_{t-2}=-\dsum%
\limits_{k=0}^{(t-9)/4}a_{4k+3},  \tag{5}
\end{equation}%
\textit{or }$t\equiv 3\func{mod}4,$

\begin{equation}
\text{ }a_{t}=-1-2\dsum\limits_{k=0}^{(t-7)/4}a_{4k+3},\text{ }%
a_{t-1}=-1-\dsum\limits_{k=1}^{(t-3)/2}a_{2k}\text{ \textit{and} \ }%
a_{t-2}=-\dsum\limits_{k=0}^{(t-7)/4}a_{4k+1};  \tag{6}
\end{equation}

\textit{(iii) \ }$\alpha ^{3}-1$\textit{\ is a unit\ }\ \textit{if and only
if} 
\begin{equation}
t\equiv 0\func{mod}3,\text{ \ }\
a_{t}=-3-2\dsum\limits_{k=1}^{(t-3)/3}a_{3k}\ \ \text{\textit{and}\ }\ \
a_{t-1}=-\dsum\limits_{\substack{ k=1 \\ k\neq 0\func{mod}3}}^{t-2}a_{k}, 
\tag{7}
\end{equation}%
\textit{or\ } 
\begin{equation}
t\equiv 1\func{mod}3,\text{ \ \ }a_{t}=-1-2\dsum%
\limits_{k=0}^{(t-4)/3}a_{3k+1}\text{ \ \ \textit{and} \ \ \ \ }%
a_{t-1}=-1-\dsum\limits_{\substack{ k=1 \\ \text{ }k\neq 1\func{mod}3}}%
^{t-2}a_{k},  \tag{8}
\end{equation}%
\textit{or} 
\begin{equation}
t\equiv 2\func{mod}3,\ a_{t}=-1-2\dsum\limits_{k=0}^{(t-5)/3}a_{3k+2}\ \ \ 
\text{\textit{and\ }\ }\ a_{t-1}=-1-\dsum\limits_{\substack{ k=1 \\ \text{ }%
k\neq 2\func{mod}3}}^{t-2}a_{k}.  \tag{9}
\end{equation}

\bigskip\ \ 

\textbf{Proof. (i) }Since $\func{Norm}(\alpha ^{2}-1)=\func{Norm}(\alpha -1)%
\func{Norm}(\alpha +1),$ $\func{Norm}(\alpha +1)=S(-1),$ and $\alpha +1$ is
a unit \ if and only if $S(-1)=1,$ we see that 
\begin{equation}
\alpha ^{2}-1\text{ is a unit}\Leftrightarrow S(1)=-1=-S(-1),  \tag{10}
\end{equation}%
and (4) follows immediately from the relation (10) and the fact that $%
(S(1)=-S(-1)=-1)\Leftrightarrow (S(1)+S(-1)=0$ and $S(1)-S(-1)=-2).$ Also,
as mentioned in the proof of [6, Proposition 3.3], we easily get from (10)
that $t\equiv 1\func{mod}2,$ whenever $\alpha ^{2}-1$ is a unit. Indeed, if $%
\alpha ^{2}-1$ is a unit, then $\alpha -1$ is a unit and $(3)$ gives that $%
a_{t}$ is odd. On the other hand, if $t$ is even then (10) yields $%
4+4a_{2}+\cdot \cdot \cdot +4a_{t-2}+2a_{t}=S(1)+S(-1)=0,$ and so $a_{t}$ is
even.

\smallskip

\textbf{(ii)} First notice that $\func{Norm}(\alpha ^{4}-1)=\func{Norm}%
(\alpha ^{2}+1)\func{Norm}(\alpha ^{2}-1),$ and $\func{Norm}(\alpha
^{2}+1)=\dprod\limits_{k=1}^{t}(\alpha _{k}^{2}+1)(\alpha
_{k}^{-2}+1)=\dprod\limits_{k=1}^{t}(\alpha _{k}+i)(\alpha
_{k}^{-1}+i)(\alpha _{k}-i)(\alpha _{k}^{-1}-i)=S(-i)S(i)=\left\vert
S(i)\right\vert ^{2},$ where $i^{2}=-1.$ Also, a simple calculation gives
that $S(i)=(2(a_{1}-a_{3}+a_{5}-\cdot \cdot \cdot -a_{t-2})+a_{t})i$\ \
(resp. $S(i)=(2(a_{1}-a_{3}+a_{5}-\cdot \cdot \cdot +a_{t-2})-a_{t})i),$
when $t\equiv 1\func{mod}4$ \ (resp. when $t\equiv 3\func{mod}4).$

Now, suppose $\alpha ^{4}-1$ unit. Then, $\ \alpha ^{2}-1$ is a unit and so
by (4) $t$ is odd and $a_{t}=-1-2(a_{1}+a_{3}+\cdot \cdot \cdot +a_{t-2}).$
Assuming $t\equiv 1\func{mod}4$ (resp. $t\equiv 3\func{mod}4),$ we obtain
from the last two equalities%
\begin{equation*}
S(i)=-(1+4(a_{3}+a_{7}+\cdot \cdot \cdot +a_{t-2}))i\ \ \ \ \ \ \text{(resp. 
}S(i)=(1+4(a_{1}+a_{5}+\cdot \cdot \cdot +a_{t-2}))i).
\end{equation*}%
Moreover, as $\alpha ^{2}+1$ is also a unit it follows from the above that 
\begin{equation*}
(1+4(a_{3}+a_{7}+\cdot \cdot \cdot +a_{t-2}))^{2}=\left\vert S(i)\right\vert
^{2}=1\text{\ \ \ (resp. }(1+4(a_{1}+a_{5}+\cdot \cdot \cdot
+a_{t-2}))^{2}=\left\vert S(i)\right\vert ^{2}=1),
\end{equation*}%
$a_{3}+a_{7}+\cdot \cdot \cdot +a_{t-2}=0$ (resp. $a_{1}+a_{5}+\cdot \cdot
\cdot +a_{t-2}=0)$ and $S(i)=-i^{t}.$ Therefore, $a_{t-2}=-a_{3}-a_{7}-\cdot
\cdot \cdot -a_{t-6}$ (resp. $a_{t-2}=-a_{1}-a_{5}-\cdot \cdot \cdot
-a_{t-6})$ and this equality together with the two ones in (4) yield the
desired relations in (5) (resp. in (6)).

Conversely, if the three coefficients $a_{t-2},$ $a_{t-1}$ and $a_{t}$ of $S$
\ satisfy the equalities in (5) (resp. in (6)), then a direct calculation
gives that $S(1)=-S(-1)=-1,$ $S(i)=-i^{t}$ and so by the above $\alpha ^{4}-1
$ is a unit.

\smallskip

\textbf{(iii)} We restrict the proof only to the case $t\equiv 0\func{mod}3,$
since similar arguments lead to the other cases. As above, we have $\func{%
Norm}(\alpha ^{3}-1)=\func{Norm}(\alpha ^{2}+\alpha +1)\func{Norm}(\alpha
-1),$ $(\alpha ^{2}+\alpha +1)=(\alpha -j)(\alpha -1/j),$ where $j=e^{i2\pi
/3},$ $S(j)=2-a_{1}-a_{2}+2a_{3}+\cdot \cdot \cdot -a_{t-2}-a_{t-1}+a_{t},$
and $\alpha ^{2}+\alpha +1$ is a unit if and only if $S(j)\overline{S(j)}=1.$

Suppose $\alpha ^{3}-1$ is a unit. Then, $\alpha -1$ is a unit and so, by
(3), $a_{t}=$ $-3-2(a_{1}+\cdot \cdot \cdot +a_{t-1}).$ It follows from the
last two equalities that $S(j)=-1-3\dsum\limits_{\substack{ k=1 \\ k\equiv
\pm 1\func{mod}3}}^{t-1}a_{k}.$ Furthermore, as $\alpha ^{2}+\alpha +1$ is a
unit we have by the above $(1+3\dsum\limits_{\substack{ k=1 \\ k\equiv \pm 1%
\func{mod}3}}^{t-1}a_{k})^{2}=\left\vert S(j)\right\vert ^{2}=1,$ $%
\dsum\limits_{\substack{ k=1 \\ k\equiv \pm 1\func{mod}3}}^{t-1}a_{k}=0$ and 
$S(j)=-1.$ Hence, $a_{t-1}=-\dsum\limits_{\substack{ k=1 \\ k\equiv \pm 1%
\func{mod}3}}^{t-2}a_{k}$ \ and so $a_{t}=-3-2\dsum%
\limits_{k=1}^{(t-3)/3}a_{3k}.$ Finally, if the last two equalities in (7)
are true, then a direct calculation shows that $S(1)=S(j)=-1,$ and hence $%
\alpha ^{3}-1$ is a unit. 
\endproof%

\bigskip

\textbf{Some examples.} \textbf{(1)} Suppose $\alpha ^{2}-1$ is a unit.
Then, (4) asserts that $t\geq 3,$ and the minimal polynomial of $\alpha $ is
of the form 
\begin{equation*}
F_{a}(x):=(x^{6}+1)-a(x^{5}+x)-(x^{4}+x^{2})-(1-2a)x^{3},
\end{equation*}%
for some $a\in \mathbb{Z},$ when $t=3.$ In fact, it was signaled in [6] \
that $F_{a}$ is a Salem polynomial for all $a\geq 0;$ hence there are
infinitely many Salem number of degree $6$ such that $\alpha ^{2}-1$ is a
unit, and this confirms Theorem 1(ii) for $(n,t)=(2,3).$

\textbf{(2)} Similarly, if $\alpha ^{4}-1$\textit{\ }is a unit, then $t$ is
an odd integer at least $3,$ and the relation (6) gives for $t=3$ that $%
(a_{1},a_{2},a_{3})=(0,-1,$ $-1),$ leading to the polynomial $F_{0}.$ Hence,
the Salem number $\alpha _{0}=1.401...$ root of $F_{0}$ is the unique Salem
number $\alpha $ of degree $6$ such that $\alpha ^{4}-1$ is a unit ($\alpha
_{0}$ is the smallest Salem number of degree $6);$ thus Theorem 1(iii) is
optimal for $n=4.$

\textbf{(3)} By the same way we obtain, using (9), that there is a unique
quartic Salem number $\alpha $ such that $\alpha ^{3}-1$ is a unit, namely $%
\alpha =1.422...$ root of $(1+x^{4})-(x+x^{3})-x^{2}$ ($\alpha $ is the
smallest Salem number of degree $4),$ and so Theorem 1(i) is sharp for $n=3.$
Notice also that each polynomial of the form 
\begin{equation*}
G_{a}(x):=(x^{6}+1)-a(x^{5}+x)+a(x^{4}+x^{2})-3x^{3}\in \mathbb{Z}[x],
\end{equation*}%
where $a\geq 3,$ is the minimal polynomial of a Salem number $\alpha $ such $%
\alpha ^{3}-1$ is a unit. Indeed, a short computation gives that the
polynomial $g_{a}(x):=(x-2)(x+1)(x-a+1)-1$ satisfies $%
g_{a}(-1)=-1<0<g_{a}(0)=2a-3,$ $g_{a}(2)=-1<0<g_{a}(1)=2a-5,$ $%
g_{a}(a-1)=-1<0<g_{a}(a)=(a-2)(a+1)-1,$ and $g_{a}$ is a cubic Salem trace
polynomial. Then, from the identity $G_{a}(x)=x^{3}g_{a}(x+1/x)$ and the
relation (7) we get the claimed result, and hence Theorem 1(i) is true for $%
n=t=3.$

\textbf{(4) }Consider the family of polynomials 
\begin{equation*}
H_{a}(x):=(x^{10}+1)-a(x^{9}+x)-a(x^{8}+x^{2})-(1-a)(x^{6}+x^{4})-(1-2a)x^{5},
\end{equation*}%
where $a\in \mathbb{Z}.$ A simple calculation shows that the polynomial $%
h_{a}(x):=$ $x(x^{2}-4)(x+1)(x-(a+1))-1$ satisfies $%
H_{a}(x)=x^{5}h_{a}(x+1/x),$ and from Theorem 3 with $%
(n,t,C_{4}(x),C(x))=(4,5,x,x+1),$ we obtain for large values of $a$ that $%
H_{a}$ is a minimal polynomial of a Salem number $\alpha $ such that $\alpha
^{4}-1.$ Incidentally, the coefficients of $H_{a}$ satisfy the equalities in
(9), and so $\alpha ^{3}-1$ is also a unit. Hence, there are infinitely many
Salem numbers of degree $10$ such that $\alpha ^{n}-1$ is a unit for all $%
n\leq 4.$

\bigskip

In terms of the Salem trace polynomial $T$ \ of $\alpha ,$ Proposition 1 may
rewritten as follows.

\medskip

\textbf{Proposition 2. }\textit{We have the following assertions. }%
\begin{equation}
\alpha -1\text{ \textit{is a unit}}\Leftrightarrow \text{\textit{\ }}T(2)=-1,
\tag{11}
\end{equation}

\begin{equation}
\alpha ^{2}-1\text{ \textit{is a unit}}\Leftrightarrow (t\equiv 1\func{mod}2%
\text{ \textit{and }}T(-2)=T(2)=-1),  \tag{12}
\end{equation}%
\begin{equation}
\alpha ^{3}-1\text{\textit{\ is a unit}}\Leftrightarrow T(-1)=T(2)=-1, 
\tag{13}
\end{equation}%
\textit{and} 
\begin{equation}
\alpha ^{4}-1\text{ \textit{is a unit}}\Leftrightarrow (t\equiv 1\func{mod}2%
\text{ \textit{and} }T(-2)=T(0)=T(2)=-1).  \tag{14}
\end{equation}

\medskip

\textbf{Proof. }Using the relation (1), it is easy to see that $%
(2)\Leftrightarrow (11),$ $(4)\Leftrightarrow (12),$ ($(7),$ $(8)$ and $(9)$)%
$\Leftrightarrow (13),$ ($(5)$ and $(6)$)$\Leftrightarrow (14),$ and so
Proposition 2 follows from Proposition 1.%
\endproof%

\bigskip\ 

A conclusion similar to the previous ones is also obtained when $n=6.$

\bigskip

\textbf{Proposition 3. }\textit{Let }$\alpha \in \mathbb{T}.$ \textit{Then,} 
\begin{equation}
\alpha ^{6}-1\text{ \textit{is a unit}}\Leftrightarrow (t\equiv 1\func{mod}2%
\text{ \textit{and} }T(-2)=T(-1)=T(1)=T(2)=-1).  \tag{15}
\end{equation}

$\bigskip $

\textbf{Proof. }The same arguments as in the proof of Proposition 1 and
Proposition 2 lead to the result.%
\endproof%

\bigskip\ 

It follows immediately from Proposition 3, when $\alpha ^{6}-1$ is a unit,
that $\deg (T)=t\geq 5,$ and so Theorem 1(ii) is also optimal for $n=6.$
Moreover, if $Q(x)$ and $R(x)$ denote respectively the quotient and the
remainder of the Euclidean division of $T(x)$ by $(x^{2}-1)(x^{2}-4),$ then $%
T(x)=(x^{2}-1)(x^{2}-4)Q(x)+R(x),$ $\deg (R)\leq 3,$ and so $R(x)=-1,$ as
(15) yields $R(-2)=R(-1)=R(1)=R(2)=-1,$ i. e., $%
T(x)=C_{6}(x)(x^{2}-4)Q(x)-1, $ where $(x^{2}-1)=C_{6}(x).$

Identically, we obtain from (11) and (13) (from (12) and (14)) that the
Salem trace polynomial of a Salem number $\alpha $ such that $\alpha ^{n}-1$
is a unit for some $n\in \{1,3\}$ (resp. $n\in \{2,4\})$ is of the form $%
C_{n}(x)(x-2)Q(x)-1$ (resp. $C_{n}(x)(x^{2}-4)Q(x)-1)$ for some monic
polynomial $Q(x)\in \mathbb{Z}[x]$ (resp. $Q(x)\in \mathbb{Z}[x]$ with odd
degree).

Let us now consider the case $n=5.$ Then, Theorem 1(i) asserts that for any
integer $t\geq 4$ there are infinitely many Salem numbers $\alpha $ of
degree $2t$ such that $\alpha ^{5}-1$ is a unit. In fact, this result
remains true for $t=3.$

\bigskip

\textbf{Proposition 4. }\textit{If }$t\leq 3$ \textit{is an integer for
which there are infinitely many Salem numbers }$\alpha $\textit{\ of degree }%
$2t$ \textit{such that }$\alpha ^{5}-1$ \textit{is a unit, then }$t=3.$

\bigskip\ 

\textbf{Proof. }First notice that if $\alpha $ is a Salem number of degree $%
2t\leq 6$ and $\alpha ^{5}-1$\textit{\ }is a unit, then $t=3.$ Indeed,
assume on the contrary that $t=2,$ and set $T(x):=x^{2}+ax+b\in \mathbb{Z}%
[x] $ for the trace polynomial of $\alpha .$ Then, the relation (1) gives $%
T(2)=-1$ and $T(2\cos (2\pi /5))T(2\cos (4\pi /5))\in \{-1,1\},$ as $\alpha
-1$ and $\alpha ^{4}+\alpha ^{3}+\alpha ^{2}+\alpha +1$ are units. Hence, $%
b=-5-2a,$ $5a^{2}+20a+11\in \{-1,1\},$ and so there is no quartic Salem
number $\alpha $ such that $\alpha ^{5}-1$\textit{\ }is a unit, since the
last relation is not true for all $a\in \mathbb{Z}.$

Now suppose $\deg (\alpha )=6$ and $\alpha ^{5}-1$ unit. Then, writing the
trace polynomial of $\alpha $ in the form $%
T(x)=(x^{2}+x-1)(x-2)+ax^{2}+bx+c, $ where $(a,b,c)\in \mathbb{Z}^{3}$ and $%
x^{2}+x-1=(x-2\cos (2\pi /5))(x-2\cos (4\pi /5))=C_{5}(x),$ we have as above 
$T(2)=-1$ and $T(2\cos (2\pi /5))T(2\cos (4\pi /5))\in \{-1,1\}.$ Hence, $%
c=-1-2b-4a,$ $5(a^{2}+b^{2}+a+b+3ab)+1\in \{-1,1\},$ and so 
\begin{equation}
a^{2}+b^{2}+a+b+3ab=0,  \tag{16}
\end{equation}%
since $5(a^{2}+b^{2}+a+b+3ab)+1\neq -1.$ On the other hand, a direct
computation gives that the resultant of a polynomial of the form $%
P_{(a,b)}(x):=(x^{2}+x-1)(x-2)+ax^{2}+bx-(1+2b+4a),$ where the integers $a$
and $b$ satisfy (16), and the polynomial $(x^{2}+x-1)(x-2)$ is equal to $-1.$

Consequently, to show Proposition 4 it is enough to prove that (16) has
infinitely many solutions $(a,b)\in \mathbb{Z}^{2}$ and the corresponding
polynomials $P_{(a,b)}$ are Salem trace polynomials. For this purpose define
a sequence $(a_{n},b_{n})_{n\geq 0}$ as follows: $a_{0}:=0,$ 
\begin{equation}
b_{n}:=\frac{-(1+3a_{n})+\sqrt{5a_{n}^{2}+2a_{n}+1}}{2},\text{ }  \tag{17}
\end{equation}%
and%
\begin{equation}
a_{n+1}:=\frac{-(1+3b_{n})-\sqrt{5b_{n}^{2}+2b_{n}+1}}{2}.  \tag{18}
\end{equation}%
Then, $(a_{0},b_{0})=(0,0),$ $(a_{1},b_{1})=(-1,2),$ $(a_{2},b_{2})=(-6,15),$
and (17) implies that $b_{n}$ is a root of $Q(x,a_{n}),$ where $%
Q(x,y)=:x^{2}+x(1+3y)+y+y^{2}=x^{2}+y^{2}+x+y+3xy\in \mathbb{Z}[x,y].$

For each $n\in \mathbb{N}$ consider the relation, say $\Re _{n},$ 
\begin{equation*}
(a_{k},b_{k})\in \mathbb{Z}^{2},\text{ \ \ }a_{k}<a_{k-1}\text{\ \ \ and \ \ 
}b_{k-1}<b_{k},\text{\ for all }k\in \{1,...n\}\text{ }
\end{equation*}%
Clearly, $\Re _{1}$ is true. Suppose that $\Re _{n}$ holds for some $n.$
Because $Q(x,y)=Q(y,x)$ and $Q(b_{n},a_{n})=0$ we see that $%
Q(a_{n},b_{n})=0, $ $a_{n}$ is a root of $Q(x,b_{n})$ and so $x-\alpha
:=Q(x,b_{n})/(x-a_{n})\in \mathbb{Z}[x].$ On the other hand, (18) asserts
that $a_{n+1}$ is a root $Q(x,b_{n})$ and so $a_{n+1}\in \{a_{n},\alpha
\}\subset \mathbb{Z},$ i. e., $Q(a_{n+1},b_{n})=0.$ Further, the
inequalities $0\leq b_{n-1}<b_{n}$ yield, by (18), that $a_{n}>$ $a_{n+1}$
and $a_{n+1}=\alpha .$ Identically, the equality $Q(a_{n+1},b_{n})=0$ yields 
$b_{n}$ is a root of $Q(x,a_{n+1}),$ $x-\beta :=Q(x,a_{n+1})/(x-b_{n})\in 
\mathbb{Z}[x],$ and $b_{n+1}\in \{b_{n},\beta \}\subset \mathbb{Z},$ as $%
b_{n+1}$ is a root of $Q(x,a_{n+1}).$ Also, the inequalities $-1\geq
a_{n}>a_{n+1}$ yield, by (17), that $b_{n}<b_{n+1},$ $b_{n+1}=\beta ,$ and
hence $\Re _{n+1}$ is true. Therefore, all pairs $(a_{n},b_{n})$ belong to $%
\mathbb{Z}^{2},$ are distinct and satisfy (16).

To conclude it is enough to verify that $P_{(a_{n},b_{n})}$ is a Salem trace
polynomial. Clearly, $P_{(a_{n},b_{n})}(0)=1-2b_{n}-4a_{n}$ and using \ (17)
it is easy to see that $P_{(a_{n},b_{n})}(0)<0$ when $a_{n}\leq -2;$ thus $%
P_{(a_{n},b_{n})}(0)<0$ for all $n\geq 2,$ because we have in this case by
the above $a_{n}\leq a_{2}=-6.$ In a similar way we get $%
P_{(a_{n},b_{n})}(1)>0,$ $P_{(a_{n},b_{n})}(2)<0,$ and so the cubic
polynomial $P_{(a_{n},b_{n})}$ has a root in each of the intervals $(0,1),$ $%
(1,2)$ and $(2,\infty ).$ Consequently, $P_{(a_{n},b_{n})}$ is a Salem trace
polynomial for all $n\geq 2,$ since the interval $(0,2)$ can not contain two
quadratic algebraic integers which are conjugates. Finally, a direct
computation shows that the polynomials $P_{(a_{0},b_{0})}$ and $%
P_{(a_{1},b_{1})}$ are also Salem trace polynomials.%
\endproof%

\bigskip

\begin{center}
\bigskip

\textbf{3. Proofs of the theorems}
\end{center}

\textbf{Proof of Theorem 1. }Let $n$ be a natural number satisfying $n\equiv
1\func{mod}2$ (resp. $n\equiv 2\func{mod}4,$ $n\equiv 0\func{mod}4).$ As
mentioned in the introduction the main tool of this proof is Theorem 3, and
so it is enough to find for each integer $d\geq 0$ a monic separable
polynomial $D(x)\in \mathbb{Z}[x]$ of degree $d$ \ (resp. of degree $2d,$ of
degree $2d+1)$ with all roots (if any) in $(-2,2)$ and such that $\gcd
(D(x),C_{n}(x))=1.$ For this purpose we shall mainly use the first kind
Chebyshev polynomials $t_{k},$ defined by the relation 
\begin{equation*}
t_{k}(2\cos \theta )=2\cos k\theta ,\text{ }\forall (k,\theta )\in \mathbb{N}%
\times \lbrack 0,\pi ].
\end{equation*}%
Clearly, $t_{1}(x)=x,$ $t_{2}(x)=x^{2}-2,$ and $%
t_{k+2}(x)=xt_{k+1}(x)-t_{k}(x),$ $\forall k\in \mathbb{N}.$ Also, $%
t_{k}(x)\in \mathbb{Z}[x],$ $\deg (t_{k})=k,$ $t_{k}$ is monic, and $t_{k}$
is separable, since its roots are $2\cos (\frac{\pi }{2k})>2\cos (\frac{\pi 
}{2k}+\frac{\pi }{k})>\cdot \cdot \cdot >2\cos (\frac{\pi }{2k}+\frac{\pi
(k-1)}{k}).$ Moreover, we have the following. 

\medskip

\textbf{Lemma 1.}\textit{\ If }$n\neq 0\func{mod}4,$\textit{\ then }$\gcd
(t_{k},C_{n})=1$\textit{\ for all} $k\in \mathbb{N}.$

\bigskip

\textbf{Proof.}\textit{\ }Let $(k,n)\in \mathbb{N}^{2}$ such that $\gcd
(t_{k},C_{n})\neq 1.$ Then, $n\geq 3,$ as $C_{1}(x)=C_{2}(x)=1,$ and there
is a root of $C_{n}$ which is also a root of $t_{k},$ i. e., there is a
natural number $l<n/2$ such that $t_{k}(2\cos (2l\pi /n))=0.$ Hence, $2\cos
(2kl\pi /n)=0$ and so $2kl\pi /n=\pi /2+\pi m$ for some $m\in \mathbb{Z};$
thus $\ 4kl=n(1+2m)$ and $n\equiv 0\func{mod}4.$%
\endproof%

\bigskip

\textbf{Lemma 2.}\textit{\ Let }$(n,m)\in \mathbb{N}^{2}.$ \textit{Then,} $%
\gcd (C_{n},C_{m})=1\Leftrightarrow \gcd (n,m)\in \{1,2\}.$

\medskip

\textbf{Proof.}\textit{\ }The result follows immediately from the definition
of the polynomials $C_{n}$ and the well known identity $\gcd
(x^{n}-1,x^{m}-1)=x^{\gcd (n,m)}-1.$%
\endproof%

\bigskip

Let us now complete the proof of Theorem 1. If we set $t_{0}(x):=1$ and $%
D(x):=t_{d}(x)$ \ (resp. and $D(x):=t_{2d}(x),$ and $D(x):=C_{4d+3}(x)),$
when $n\equiv 1\func{mod}2$ (resp. when $n\equiv 2\func{mod}4,$ when $%
n=2^{s} $\textit{\ }for some integer\textit{\ }$s\geq 2),$ then we easily
get by Lemma 1 (resp. by Lemma 1, by Lemma 2) that Theorem 1(i) (resp.
Theorem 1(ii), Theorem 1(iii)) is true.

To show the last assertion in Theorem 1, suppose $n\equiv 4\func{mod}8,$ $%
n\neq 0\func{mod}3$ and\ $D(x):=(x-1)t_{2d}(x).$ Then, $D(x)\in \mathbb{Z}%
[x],$ and $D$ is monic and separable, since $C_{6}(x)=(x-1)(x+1)$ and Lemma
1 yields $\gcd (t_{k},C_{6})=1$ for all $k.$ To conclude notice that $\gcd
(C_{n},D)=1,$ i. e., $1\neq 2\cos (2l\pi /n)$ and $t_{2d}(2\cos (2l\pi
/n))\neq 0$ for all $l\in \{1,...,(n-2)/2\}.$ Indeed, if $1=2\cos (2l\pi /n)$
for some $l\in \{1,...,(n-2)/2\},$ then $\cos (2l\pi /n)=\cos (\pi /3),$ $%
2l\pi /n=\pi /3,$ $6l=n$ and so $n\equiv 0\func{mod}3.$ Similarly, if $%
t_{2d}(2\cos (2l\pi /n))=0$ and $l\in \{1,...,(n-2)/2\},$ then $\cos (4dl\pi
/n)=0$ and $4dl\pi /n=\pi /2+\pi m$ \ \ for some $m\in \mathbb{Z};$ thus $%
8dl=n(1+2m)$ and so $n\equiv 0\func{mod}8.$%
\endproof%

\ \ 

\bigskip

\bigskip

\textbf{Proof of Theorem 2.} Let $n\in \mathbb{N}.$ Without loss of
generality we may restrict the proof to the case $n\equiv 0\func{mod}4,$
since for the other values of $n$ Theorem 2 is an immediate corollary of
Theorem 1. Then, there is an infinite subset $V$ of $\mathbb{N}$ such that $%
\gcd (n,4v+3)=1$ for each $v\in V.$ This assertion follows (for instance)
from Dirichlet's Theorem, saying that there are infinitely many prime
numbers of the form $4v+3.$ Then, Lemma 2 \ gives that $\gcd
(C_{n},C_{4v+3})=1,$ for each $v\in V,$ \ and it follows by Theorem 3 with $%
D:=C_{4v+3}$ that there are infinitely many Salem numbers $\alpha $ of
degree $2t=2(2v+1+(n+4)/2)=4v+n+6$ such that $\alpha ^{n}-1$ is a unit.
Setting $I:=\{2v+3+n/2\mid v\in V\}$ we obtain the desired result. 
\endproof%

\ \ 

\medskip

\textbf{Proof of Theorem 3. }Set $P(x):=C_{n}(x)D(x)(x-2)(x-a)$ \ \ (resp. $%
P(x):=C_{n}(x)D(x)(x^{2}-4)(x-a))$\ when $n$ is odd (resp. when $n$ is even)
and 
\begin{equation}
R_{a}(x)=R(x):=P(x)-1.  \tag{19}
\end{equation}%
Clearly, the polynomial product $C_{n}D$ is separable with all zeros in the
interval $(-2,2),$ since so are the coprime polynomials $C_{n}$ and $D.$
Also, $\deg (D)=t-(n+3)/2\geq 0,$ $\deg (C_{n})=$ $(n-1)/2$ (resp. $\deg
(D)=t-(n+4)/2)\geq 0,$ $\deg (C_{n})=$ $(n-2)/2),$ and $\deg (P)=\deg (R)=$ $%
t\geq 2$ (resp. $t\geq 3).$

Suppose $a\geq 3.$ Then, the polynomial $P$ is separable and its roots, say $%
\beta _{1},\beta _{2},...,$ $\beta _{t},$\ may be labelled so that 
\begin{equation*}
-2\leq \beta _{1}<\cdot \cdot \cdot <\beta _{t-1}=2<\beta _{t}=a.\text{ \ \ }
\end{equation*}%
\ To make clear the remaining part of the proof consider the following two
cases.

\ 

\textbf{Case }$t$\textbf{\ odd. }Then $t\geq 3.$ For each $k\in
\{1,...,(t-1)/2\}$ fix an element $\gamma _{k}$ of the interval $(\beta
_{2k-1},\beta _{2k}).$ Then, $(\gamma _{k}-\beta _{1})>0,$ $...,$ $(\gamma
_{k}-\beta _{2k-1})>0,$ $(\gamma _{k}-\beta _{2k})(\gamma _{k}-\beta
_{2k+1})\cdot \cdot \cdot (\gamma _{k}-\beta _{t})>0,$ and so $P(\gamma
_{k})>0.$ Moreover, if $a>A_{k},$ where 
\begin{equation}
A_{k}:=\left\vert \gamma _{k}\right\vert +\frac{1}{\left\vert (\gamma
_{k}-\beta _{1})(\gamma _{k}-\beta _{2})\cdot \cdot \cdot (\gamma _{k}-\beta
_{t-1})\right\vert },  \tag{20}
\end{equation}%
then $\left\vert \gamma _{k}-\beta _{t}\right\vert \geq $ $\beta
_{t}-\left\vert \gamma _{k}\right\vert =a-\left\vert \gamma _{k}\right\vert
>1/\left\vert (\gamma _{k}-\beta _{1})(\gamma _{k}-\beta _{2})\cdot \cdot
\cdot (\gamma _{k}-\beta _{t-1})\right\vert ,$ and hence $P(\gamma _{k})>1.$
It follows, by (19), that 
\begin{equation*}
R(\gamma _{k})>0>-1=R(\beta _{2k-1})=R(\beta _{2k})
\end{equation*}%
and the polynomial $R$ has two roots, say $\eta _{2k-1}$ and $\eta _{2k}$
such that 
\begin{equation*}
\beta _{2k-1}<\eta _{2k-1}<\gamma _{k}<\eta _{2k}<\beta _{2k}.
\end{equation*}%
Therefore, for each integer $a>A:=\max \{A_{1},...,A_{(t-1)/2}\}$ the
polynomial $R$ has $(t-1)$ distinct roots belonging to the interval $(-2,2)$
and the remaining root, say $\eta _{a,t},$ belongs to the interval $%
(a,\infty ),$ as the monic polynomial $R$ satisfies $R(a)=-1<0.$ In fact,
since the monic polynomial $C_{n}D$ has no root greater than $2,$ we have $%
C_{n}(a+1)D(a+1)\geq 1,$ $R(a+1)>0$ and so $\eta _{a,t}\in (a,a+1).$

To complete the proof of this case, suppose $a>A$ and set $%
R_{a}(x)=Q_{a}(x)M_{a}(x),$ where $M_{a}$ is the minimal polynomial of $\eta
_{a,t}.$ Then, $Q_{a}(x)\in \mathbb{Z}[x],$ $\deg (Q_{a})\leq t-2$ (recall
that $\eta _{a,t}\in (a,a+1)$ and so $\deg (M_{a})\geq 2),$ $Q_{a}$ is
monic, and the roots of $Q_{a}$ (when $\deg (Q_{a})\geq 1)$ belong to the
interval $(-2,2).$

Because there is a finite number of polynomials in $\mathbb{Z}[x]$ with
given degree and bounded coefficients, the number, say $q,$ of nonconstant
monic polynomials with degree at most $t-2$ and roots in $(-2,2)$ is finite.
Notice also that if $\eta $ is a common zero of $R_{a}$ and $R_{a^{\prime
}}, $ for some $(a,a^{\prime })\in (A,\infty )^{2},$ then $C_{n}(\eta
)D(\eta )(\eta -2)(\eta -a)=1=C_{n}(\eta )D(\eta )(\eta -2)(\eta -a^{\prime
})$ (resp. $C_{n}(\eta )D(\eta )(\eta ^{2}-4)(\eta -a)=1=C_{n}(\eta )D(\eta
)(\eta ^{2}-4)(\eta -a^{\prime })),$ $a=a^{\prime }$ and so $\gcd
(R_{a},R_{a^{\prime }})=1$ when $a\neq a^{\prime }.$ Consequently, there is
at most $q$ polynomials $R_{a}$ which are reducible, when $a$ runs through
the interval $(A,\infty )$ and so there is a constant $B$ such that if $a>B,$
then $R_{a}$ is irreducible, i. e., $R_{a}$ is a trace polynomial of a Salem
polynomial $S_{\alpha },$ where $\alpha +1/\alpha =\eta _{a,t}.$

Finally, notice that if $\zeta $ is a zero of $x^{n}-1,$ then $P(\zeta
+1/\zeta )=0,$ $R_{a}(\zeta +1/\zeta )=-1,$ and so, by (1), $S_{\alpha
}(\zeta )=-\zeta ^{t}.$ Hence, the absolute value of the resultant $%
\dprod\limits_{\zeta ^{n}=1}S_{\alpha }(\zeta )$ of the polynomials $x^{n}-1$
and $S_{\alpha }(x)$ is equal to $1,$ and consequently $\alpha ^{n}-1$ is a
unit.

\textbf{Case }$t$\textbf{\ even. } Since $\deg (P)=t$ \ and $P$ has no root
less than $-2,$ we have $P(-2)\geq 0.$ Also, [6, Proposition 3.3] gives that 
$n$ \ is odd. Thus, $P(-2)>0,$ $P(-2)\geq 1,$ $R(-2)\geq 0$ and $R$ has a
root $\eta _{1}\in $ $[-2,\beta _{1}),$ as $R(\beta _{1})=-1<0.$

Now, if fix for each $k\in \{1,...,(t-2)/2\},$ where $t\geq 4,$ an element $%
\gamma _{k}$ of the interval $(\beta _{2k},\beta _{2k+1}),$ then we get as
above $(\gamma _{k}-\beta _{1})>0,$ $...,$ $(\gamma _{k}-\beta _{2k})>0,$ $%
(\gamma _{k}-\beta _{2k+1})\cdot \cdot \cdot (\gamma _{k}-\beta _{t})>0,$
and hence $P(\gamma _{k})>0.$ Also, if we choose $a>A_{k},$ where $A_{k}$ is
defined by (20), then we obtain $P(\gamma _{k})>1,$

\begin{equation*}
R(\gamma _{k})>0>-1=R(\beta _{2k})=R(\beta _{2k+1})
\end{equation*}%
and so the polynomials $R$ has two roots, say $\eta _{2k}$ and $\eta _{2k+1}$
such that 
\begin{equation*}
\beta _{2k}<\eta _{2k}<\gamma _{k}<\eta _{2k+1}<\beta _{2k+1}.
\end{equation*}%
From this point we conclude in the same way as in the case where $t$ is odd.%
\endproof%

\bigskip

\begin{center}
\textbf{References}
\end{center}

[1] C. Batut, D. Bernardi, H. Cohen and M. Olivier, \textit{User's Guide to
PARI-GP}, Version 2.5.1 (2012).

[2] M. J. Bertin, A. Decomps-Guilloux, M. Grandet-Hugo, M.
Pathiaux-Delefosse and J. P. Schreiber, \textit{Pisot and Salem numbers},
Birkh\"{a}user Verlag Basel, 1992.

[3] D. W. Boyd, \textit{Small Salem numbers}, Duke Math. J. \textbf{44}
(1977), 315-328.

[4] D. W. Boyd, \textit{Pisot and Salem numbers in intervals of the real line%
}, Math. Comp. \textbf{32} (1978), 1244-1260.

[5] J. H. Evertse, \textit{On equations in S-units and the Thue-Mahler
equation}, Invent. Math. \textbf{75} (1984), 561-584.

[6] B. H. Gross and C. T. McMullen, \textit{Automorphisms of even unimodular
lattices and unramified Salem numbers}, J. Algebra. \textbf{257} (2008),
265-290.

[7] J. H. Silverman, \textit{Exceptional units and numbers of small Mahler
measure,} Exp. Math. \textbf{4} (1995), 69-83.

[8] T. Za\"{\i}mi, \textit{Remarks on certain Salem numbers,} Arab J. Math.
Sc., \textbf{7} (2001), 1- 10.

[9] T. Za\"{\i}mi, \textit{An arithmetical property of \ powers of Salem
numbers,} J. Number Theory \textbf{120} (2006), 179- 191.

[10] T. Za\"{\i}mi, \textit{Comments on Salem polynomials,} Arch. Math. 
\textbf{117} (2021), 41--51.

\bigskip

Department of Mathematics and Statistics, College of Science

Imam Mohammad Ibn Saud Islamic University (IMSIU)

P. O. Box 90950

Riyadh 11623 Saudi Arabia

Email: tmzaemi@imamu.edu.sa\textit{\ }

\bigskip

\end{document}